\documentclass{amsart}
\usepackage{amssymb} 

\usepackage{color} 

\newtheorem{corollary}{Corollary}
\newtheorem{theorem}{Theorem}
\newtheorem{lemma}{Lemma}

\begin{document}
\title{Numerical invariants of identities of unital algebras}

\author[D. Repov\v s, and M. Zaicev]
{Du\v san Repov\v s and Mikhail Zaicev}

\address{Du\v{s}an Repov\v{s} \\Faculty of Education, and
Faculty of Mathematics and Physics, University of Ljubljana,
Kardeljeva plo\v s\v cad 16, Ljubljana, 1000, Slovenia}
\email{dusan.repovs@guest.arnes.si}

\address{Mikhail Zaicev \\Department of Algebra\\ Faculty of Mathematics and
Mechanics\\  Moscow State University \\ Moscow, 119992, Russia}
\email{zaicevmv@mail.ru}

\thanks{The first author was supported by
SRA grants P1-0292-0101, J1-4144-0101 and J1-5435-0101.
The second author was supported by 
RFBR grant No 13-01-00234a}



\begin{abstract}
We study polynomial identities of algebras with adjoined external unit. For a wide class
of algebras we prove that adjoining external unit element leads to increasing
of PI-exponent precisely to 1. We also show that any real number from the
interval [2,3] can be realized as PI-exponent of some unital algebra.
\end{abstract}

\date{\today}

\maketitle

\section{Introduction}

We study numerical characteristics of polynomial identities of algebras over a field
$F$ of characteristic zero. Given an algebra $A$ over $F$, one can associate to it the
sequence $\{c_n(A)\}$ of non-negative integers called the {\it sequence of codimensions}.
If the growth of $\{c_n(A)\}$ is exponential then the limiting ration of consecutive terms 
is called PI-{\it exponent} of $A$ and written  $exp(A)$. 
In the present paper we are mostly interested what happens with PI-exponent if
we adjoin to $A$ an external unit element.

The first results in this area were proved for associative algebras. It is known that $exp(A)$
is an integer in the associative case \cite{GZ98}, \cite{GZ99}.  It was
shown in \cite{GZJA} that it follows from the proofs in \cite{GZ98}, \cite{GZ99} that either $exp(A^\sharp)
=exp(A)$ or $exp(A^\sharp)=exp(A)+1$ and both options can be realized. Here $A^\sharp$ is the algebra
$A$ with adjoined external unit.

The next result was published in \cite{ZAL},  following an example of A 5-dimensional
 algebra $A$ with $exp(A)<2$  constructed in \cite{GMZCA2009}. The point
is that in the associative or Lie case PI-exponent cannot be less than 2 (\cite{K}, \cite{M}).
For a finite dimensional Lie superalgebra, Jordan and alternative algebra PI-exponent 
is also at least 2. Starting from
the example $A$ from \cite{GMZCA2009} it was shown in \cite{ZAL}  that
$exp(A^\sharp)=exp(A)+1$. In \cite{ZAL} also the following problem
was stated: is it true that always 
either $exp(A^\sharp)=exp(A)$ or $exp(A^\sharp)=exp(A)+1$?

An example of 4-dimensional simple algebra $A $ with a fractional 
PI-exponent was constructed in \cite{BBZF}. It was also shown that
$exp(A^\sharp)=exp(A)+1$. This 
result was announced in \cite{BBZU}. It was also shown in \cite{BBZU} that if
$A$ is itself a unital algebra then $exp(A^\sharp)=exp(A)$.

In the present paper (see Theorem \ref {t1}) we shall prove that for a previously known 
series  of algebras $A_\alpha$ with $exp(A_\alpha)=\alpha$, $\alpha\in \mathbb R,
1<\alpha<2$ (see \cite{GMZAdvM}) any extended algebra $A^\sharp_\alpha$ has exponent
$\alpha+1$. That is, we shall show that there exist infinitely many algebras $A$
such that $exp(A^\sharp)=exp(A)+1$.

Another important question is the following: which real numbers can be realized as
PI-exponents of some algebra? For example, if $A$ is any associative PI-algebra
or a finite dimensional Lie or Jordan algebra then $exp(A)$ is an integer (see
\cite{GZ98}, \cite{GZ99}, \cite{Z}, \cite{GSZ}).

For unital algebras it is only known  that if $\dim A<\infty$ then $exp(A)$
cannot be less than 2. As a consequence of the main result of our paper (see
Corollary \ref{c1}) we shall obtain that for any real $\alpha\in [2,3]$ there exists
a unital algebra $B_\alpha$ such that $exp(B_\alpha)=\alpha$.

\section{Preliminaries}

Let $A$ be an algebra over a field $F$ of characteristic zero and let
$F\{X\}$ be absolutely free algebra over $F$ with a countable set of generators
$X=\{x_1,x_2,\ldots\}$. Recall that a polynomial $f=f(x_1,\ldots,x_n)$ is said to be an
{\it identity} of $A$ if $f(a_1\ldots, a_n)=0$ for all $a_1,\ldots,a_n\in A$. The set $Id(A)$  
of all polynomial identities of $A$ forms an ideal of $F\{X\}$.

Denote by $P_n$ the subspace of all multilinear polynomials in $F\{X\}$ on $x_1,\ldots, x_n$. Then
the intersection $Id(A)\cap P_n$ is the space of all multilinear identities of $A$
of degree $n$. 

Denote
$$
P_n(A)=\frac{P_n}{Id(A)\cap P_n}.
$$
A non-negative integer
$$
c_n(A)=\dim P_n(A)
$$is called the $n$th {\it codimension} of $A$. Asymptotic behavior of the sequence $\{c_n(A)\},
n=1,2,\ldots,$ is an important numerical invariant of identities of $A$. We refer readers to 
\cite{GZbook} for an account of basic notions of the theory of codimensions of PI-algebras.

If the sequence $\{c_n(A)\}$ is exponentially bounded, i.e. $c_n(A)\le a^n$ for all $n$ and for some
number $a$ (for example in case $\dim A<\infty$ and in many other cases), we can define the 
lower and the upper PI-exponents of $A$ by
$$
\underline{exp}(A)=\liminf_{n\to\infty}\sqrt[n]{c_n(A)},\quad
\overline{exp}(A)=\limsup_{n\to\infty}\sqrt[n]{c_n(A)},
$$
and (the ordinary) PI-exponent
$$
exp(A)=\lim_{n\to\infty}\sqrt[n]{c_n(A)} 
$$
provided that $\underline{exp}(A)=\overline{exp}(A)$.

In order to compute the values of codimensions we can consider symmetric group action 
on $P_n$ defined by
$$
\sigma f(x_1,\ldots,x_n)=f(x_{\sigma(1)},\ldots,x_{\sigma(n)}\quad \forall \sigma\in S_n.
$$
The subspace $P_n\cap Id(A)$ is invariant under this action and we can study the structure
of $P_n(A)$ as an $S_n$-module. Denote by $\chi_n(A)$ the $S_n$-character of $P_n(A)$
called the $n$th {\it cocharacter} of $A$. Since ${\rm char}~ F=0$  and
any $S_n$-representation is completely 
reducible, the $n$th cocharacter has the decomposition
\begin{equation}\label{e1}
\chi_n(A)=\sum_{\lambda\vdash n}m_\lambda \chi_\lambda
\end{equation}
where $\chi_\lambda$ is the irreducible $S_n$-character corresponding to the partition 
$\lambda\vdash n$ and non-negative integer $m_\lambda$ is the multiplicity of $\chi_\lambda$
in $\chi_n(A)$.

Obviously, it follows from (\ref{e1}) that
$$
c_n(A)=\sum_{\lambda\vdash n}m_\lambda \deg\chi_\lambda.
$$

Another important numerical characteristic is the $n$th {\it colength} of $A$ defined by
$$
l_n(A)=\sum_{\lambda\vdash n}m_\lambda 
$$
with $m_\lambda$ taken from (\ref{e1}). In particular, if the sequence $\{l_n(A)\}$ is polynomially bounded
as a function of $n$ while some of $\deg\chi_\lambda$ with $m_\lambda\ne 0$ are exponentially large, 
the principal part of the asymptotic of $\{c_n(A)\}$ is defined by the largest value of $\deg\chi_\lambda$ 
with non-zero multiplicity.

For studying asymptotic of codimensions it is convenient to use 
the following functions. Let
$0\le x_1,\ldots,x_d\le 1$ be real numbers such that $x_1+\cdots+x_d=1$. Denote
$$
\Phi(x_1,\ldots,x_d)=\frac{1}{x_1^{x_1}\cdots x_d^{x_d}}.
$$
If $d=2$ then instead of $\Phi(x_1,x_2)$ we will write
$$
\Phi_0(x)=\frac{1}{x^x(1-x)^{1-x}}.
$$
We assume that some of $x_1,\ldots,x_d$ can have zero values. In this case we 
assume that $0^0=1$.

Given $\lambda=(l_1,\ldots,\lambda_d)\vdash n$ we define
\begin{equation}\label{e3}
\Phi(\lambda)=\frac{1}{(\frac{\lambda_1}{n})^{\frac{\lambda_1}{n}}\cdots (\frac{\lambda_d}{n})^{\frac{\lambda_d}{n}}}.
\end{equation}
For partitions $\lambda=(\lambda_1,\ldots,\lambda_k)\vdash n$ with $k<d$ we also consider $\Phi(\lambda)$
as in (\ref{e3}), assuming $\lambda_{k+1}=\cdots=\lambda_d=0$.

The relationship  between $\deg\chi_\lambda$ and $\Phi(\lambda)$ is given by 
the following lemma.
\begin{lemma}\label{l1}{\rm (see \cite[Lemma 1]{GZJLMS})}
Let $\lambda=(\lambda_1,\ldots,\lambda_k)\vdash n$ be a partition of $n$. If $k\le d$ and $n\ge 100$ then
$$
\frac{\Phi(\lambda)^n}{n^{d^2+d}} \le \deg \chi_\lambda \le n\Phi(\lambda)^n.
$$
\end{lemma}
\hfill $\Box$

Now we investigate how the value of $\Phi(x_1,\ldots, x_d)$ increases  after adding one extra variable.
\begin{lemma}\label{l3}
Let
$$
\Phi(x_1,\ldots,x_d)=\frac{1}{x_1^{x_1}\cdots x_d^{x_d}},
\quad 0\le x_1,\ldots,x_d,\quad x_1+\cdots+x_d=1,
$$
and let $\Phi(z_1,\ldots,z_d)=a$ for some fixed $z_1,\ldots,z_d$. Then
$$
\max_{{{0\le t\le 1}}}\{\Phi(y_1,\ldots y_d,1-t)\vert y_1=tz_1,\ldots, y_d=tz_d\}=a+1.
$$
Moreover, the maximal value is achieved if $t=\frac{a}{a+1}$.
\end{lemma}
{\em Proof.} Consider
$$
g(t)=\ln \Phi^{-1}(tz_1,\ldots, tz_d, 1-t).
$$
Then
$$
g(t)=t\ln t+(1-t)\ln(1-t)-t\ln a.
$$
Hence its derivative is equal to
$$
g'(t)=\ln\frac{t}{(1-t)a}
$$
and $g'(t)=0$ if and only if $t=(1-t)a$, that is $t=\frac{a}{a+1}$. It is not difficult to check that $g$
has the minimum at this point.

Now we compute the value of $g$:
$$
g(\frac{a}{a+1})=\frac{a}{a+1}\ln \frac{a}{a+1}+\frac{1}{a+1}\ln \frac{1}{a+1}-\frac{a}{a+1}\ln a=\ln B,
$$
where
$$
B=(\frac{a}{a+1})^{\frac{a}{a+1}}(\frac{1}{a+1})^{\frac{1}{a+1}} a^{-\frac{a}{a+1}}=\frac{1}{a+1}.
$$
Hence $\Phi_{max}=B^{-1}=a+1$, and we have completed the proof.
\hfill $\Box$

Next lemma shows what happens with $\Phi(\lambda)$ when we insert an extra row in Young diagram $D_\lambda$.

\begin{lemma}\label{l4}
Let $\gamma$ be a positive real number and let $\lambda=(\lambda_1,\ldots,\lambda_d)$ be a partition of $n$ such 
that $\frac{\lambda_1}{n},\ldots, \frac{\lambda_d}{n}\ge\gamma$. Then for any $\varepsilon>0$ there
exist $n'=kn$ and a partition $\mu\vdash n'$, $\mu=(\mu_1,\ldots,\mu_{d+1})$ such that for some integers
$1\le i\le d+1$ and $q\ge 1$ the following conditions hold:
\begin{itemize}
\item[1)]
$\mu_j=q\lambda_j$ for all $j\le i-1$;
\item[2)]
$\mu_{j+1}=q\lambda_j$ for all $j\ge i$; and
\item[3)]
 $\vert \Phi(\lambda)-\Phi(\mu)+1\vert <\varepsilon$.
\end{itemize}
Moreover, $k$ does not depend on $\lambda$ and $n$.
\end{lemma}
{\em Proof.} Denote
$$
z_1=\frac{\lambda_1}{n},\ldots,z_d=\frac{\lambda_d}{n}
$$
and $a=\Phi(z_1,\ldots,z_d)=\Phi(\lambda)$. By Lemma \ref{l3}
\begin{equation}\label{e3a}
\Phi(tz_1,\ldots,tz_d, 1-t)=a+1
\end{equation}
if $t=\frac{a}{a+1}$. It is not difficult to check that $1\le \Phi(x_1,\ldots,x_d)\le d$, hence 
$\frac{1}{d+1}\le 1-t\le\frac{1}{2}$. 

Note that $\Phi=\Phi(x_1,\ldots,x_{d+1})$ can be viewed as a function of $d$ independent indeterminates
$x_1,\ldots,x_d$. Conditions $0<\gamma\le x_1,\ldots, x_d$ and $\frac{1}{d+1}\le x_{d+1}
\le \frac{1}{2}$ define a compact domain $Q$ in 
$\mathbb{R}^d$ since $x_{d+1}=1-x_1-\cdots- x_d$. Since $\Phi$ is continuous on $Q$ there exists an integer $k$
such that
$$
\vert \Phi(x_1,\ldots,x_d,x_{d+1} -  \Phi(x_1',\ldots,x_d',x_{d+1}'\vert < \varepsilon
$$
as soon as $\vert x_i-x_i'\vert < \frac{1}{k}$ for all $i=1,\ldots,d$. Clearly, $k$ does not depend on $n$ and $\lambda$. Then there exists a rational number 
$t_0=\frac{q}{k}<1$ such that $\vert t-t_0\vert < \frac{1}{k}$ and
\begin{equation}\label{e4}
\vert \Phi(t_0z_1,\ldots, t_0z_d, 1-t_0)-a-1\vert <\varepsilon.
\end{equation}
Denote $y_0=1-t_0$. Then $t_0z_i\le 1-t_0=y_0  \le t_0z_{i-1}$ for  some $i$ 
(or $y_0>t_0z_1$, or $y_0< t_0 z_d$).

Now we set $n'=kn$,
$$
\mu_1=q\lambda_1,\ldots,\mu_{i-1}=q\lambda_{i-1},
$$
$$
\mu_{i+1}=q\lambda_i,\ldots, \mu_{d+1}=q\lambda_d
$$
and $\mu_i=n(k-q)$. Then $\mu=(\mu_1,\ldots,\mu_{d+1})$ is a partition of $n'$ and
$$
\Phi(\mu)=\Phi(t_0z_1,\ldots,t_0z_d,1-t_0).
$$
In particular, $\vert \Phi(\lambda)-\Phi(\mu)-1\vert <\varepsilon$ by
(\ref{e3a}) and (\ref{e4}), and we have completed the proof of our lemma.
\hfill $\Box$

\section{Algebras of infinite words}

In this section we recall some constructions and algebras from \cite{GMZAdvM}  and their properties.
These algebras will be used for constructing unital algebras.

Let $K=(k_1,k_2,\ldots)$ be an infinite sequence of integers $k_i\ge 2$. Then the algebra $A(K)$ is defined by its basis
\begin{equation}\label{e01}
\{a,b\}\cup Z_1\cup Z_2\cup\ldots
\end{equation}
where
\begin{equation}\label{e02}
Z_i=\{z^{(i)}_j\vert 1\le  j\le k_i, \,i=1,2,\ldots\}
\end{equation}
with the multiplication table
\begin{equation}\label{e03}
z^{(i)}_1a=z^{(i)}_2,\ldots, z^{(i)}_{k_i-1}a=z^{(i)}_{k_i}, z^{(i)}_{k_i}b=z^{(i+1)}_1
\end{equation}
for all $i=1,2,\ldots$. All remaining products are assumed to be zero.

It is easy to verify (see also \cite{GMZAdvM}) that $A$ satisfies the identity $x_1(x_2x_3)=0$ and if $m_\lambda\ne 0$ in (\ref{t1})
then $\lambda=(\lambda_1)$ or $\lambda=(\lambda_1,\lambda_2)$ or $\lambda=(\lambda_1,\lambda_2,1)$. Denote by 
$W_n^{(d)}$, $d\le  n$, the subspace of the free algebra $F\{X\}$ of all homogeneous polynomials of
degree $n$ on $x_1,\ldots, x_d$. Given a PI-algebra $A$, we define
$$
W_n^{(d)}(A)=\frac{W_n^{(d)}}{W_n^{(d)}\cap Id(A)}.
$$

Recall that the height $h(\lambda)$ of a partition $\lambda=(\lambda_1,\ldots,\lambda_d)$ is equal to $d$.
We will use the following result from \cite{GMZAdvM}.
\begin{lemma}\label{l5} {\rm (\cite[Lemma 4.1]{GMZAdvM})}
Let $A$ be a PI-algebra with $n$th cocharacter $\chi_n(A)=\sum_{\lambda\vdash n} m_\lambda\chi_\lambda$.
Then for every $\lambda\vdash n$ with $h(\lambda)\le d$ we have that $m_\lambda\le \dim W_n^{(d)}(A)$.
\end{lemma}
\hfill $\Box$

Now let $w=w_1w_2\ldots$ be an infinite word in the alphabet $\{0,1\}$. Given an integer $m\ge 2$, let 
$K_{m,w}=\{k_i\}, i=1,2,\ldots$, be the sequence defined by
\begin{equation}\label{e04}
k_i=\left\{
\begin{array}{rcl}

       m,\, if\, w_i=0\\
       m+1,\, if\, w_i=1
 \end{array}
\right.
\end{equation}
and write $A(m,w)=A(K_{m,w})$.

Recall that the complexity $Comp_w(n)$ of an infinite word $w$ is the number of distinct subwords of $w$ 
of the length $n$ (see \cite{Lotair}, Chapter 1). Slightly modifying the proof of Lemma 4.2 from \cite{GMZAdvM}
we obtain:
\begin{lemma}\label{l6}
For any $m\ge 2$ and for any infinite word $w$ the following inequalities hold
$$
\dim W_n^{(d)}(A(m,w))\le d(m+1)n Comp_w(n)
$$
and
$$
l_n(A(m,w))\le n^3\dim W_n^{(3)}(A(m,w)).
$$
\end{lemma}
\hfill $\Box$

Now we fix the algebra $A(m,w)$ by choosing the word $w$. Obviously, $Comp_w(n)\le T$ for any infinite periodic
word with period $T$. It is known (see \cite{Lotair}) that $Comp_w(n)\ge n+1$ for any aperiodic word $w$.
In case $Comp_w(n)=n+1$ for all $n\ge 1$ the word $w$ is said to be {\it Sturmian}. It is also known that for any Sturmian 
or periodic word the limit
$$
\pi(w)=\lim_{n\to\infty}\frac{w_1+\cdots+w_n}{n}>0
$$ 
always exists (we always assume that a periodic word is non-zero). This limit $\pi(w)$ is called the {\it  slope} of $w$.
For any real number $\alpha\in(0,1)$ there exists a word $w$ with $\pi(w)=\alpha$ and $w$ is Sturmian or
periodic depending on whether  $\alpha$ is irrational or rational, respectively. Moreover,
$$
exp(A(m,w))=\Phi_0(\beta)=\frac{1}{\beta^\beta(1-\beta)^{1-\beta}}
$$
 for Sturmian or periodic word $w$, where $\beta=\frac{1}{m+\alpha}, \alpha=\pi(w)$ (\cite{GMZAdvM}, Theorem 5.1).
As a consequence, for any real $1\le\alpha \le 2$ there exists an algebra $A$ such that $exp(A)=\alpha$.

Finally, for any periodic word $w$ and for any $m\ge 2$ there exists a finite dimensional algebra $B(m,w)$ 
satisfying the same identities as $A(m,w)$. In particular, for any rational $0<\beta\le\frac{1}{2}$ there exists
a finite dimensional algebra $B$ with
$$
exp(B)=\Phi_0(\beta)=\frac{1}{\beta^\beta(1-\beta)^{1-\beta}}.
$$

\section{Algebra with adjoined unit}

We fix an infinite or periodic word $w$ and $m\ge 2$ and consider the algebra $A=A(m,w)$. Denote by $A^\sharp$ 
the algebra obtained from $A$ by adjoining external unit element $1$. Our main goal is to prove that $exp(A^\sharp)$ exists and that 
$$
exp(A^\sharp)=exp(A)+1.
$$
First we find a polynomial upper bound for the colength of $A^\sharp$. We start with the remark
concerning an arbitrary algebra $B$. Recall that, given an algebra $B$, $W_n^{(d)}(B)$ is the dimension of 
the space of homogeneous polynomials on $x_1,\ldots,x_d$ of total degree $n$ modulo ideal 
$Id(B)$.

\begin{lemma}\label{l7}
Let $B$ be an arbitrary algebra. Suppose that $\dim W_n^{(d)}(B)\le\alpha n^T$ for some natural $T$, 
$\alpha\in\mathbb{R}$ and for all $n\ge 1$. Then
$$
\dim W_n^{(d)}(B^\sharp)\le\alpha (n+1)^{T+d+1}.
$$
\end{lemma}
{\em Proof}. Denote by $F\{X\}^\sharp$ absolutely free algebra generated by $X$ with adjoined unit element.
First note that a multihomogeneous polynomial $f(x_1,\ldots,x_d)$ is an identity of $B^\sharp$ if all
multihomogeneous on $x_1,\ldots,x_d$ components of $f(1+x_1,\ldots,1+x_d)$ are identities of $B$.

Clearly, the number of multihomogeneous polynomials on $x_1,\ldots,x_d$ of total degree 
$k$, linearly independent modulo $Id(B)$, 
does not exceed $\dim W_k^{(d)}(B)$. On the other hand,
the number of multihomogeneous components of total degree $k$  in a free algebra
$F\{x_1,\ldots,x_d\}$ does not exceed $(k+1)^d$. Take now
$$
N=(k+1)^d\sum_{k=0}^n \dim W_k^{(d)}(B) +1
$$
assuming that $\dim W_0^{(d)}(B)=1$. Clearly,
$$
N\le 1+(n+1)^d\alpha\sum_{k=0}^n k^T <\alpha(n+1)^{T+d+1}.
$$

Given  homogeneous polynomials $f_1,\ldots, f_{N+1}$ on $x_1,\ldots,x_d$
of degree $n$, consider their linear combination $f=\lambda_1 f_1+ 
\cdots+ \lambda_{N+1} f_{N+1}$ with unknown coefficients 
$\lambda_1,\cdots, \lambda_{N+1}$. The assumption that some multihomogeneous
component of $f(1+x_1,\ldots, 1+x_d)$ is an identity of $B^\sharp$ is equivalent to
some linear equation on $\lambda_1,\cdots, \lambda_{N+1}$. Hence  the
condition that all multihomogeneous components of $f(1+x_1,\ldots, 1+x_d)$
are identities of $B$
leads to at most $N$ linear equations on $\lambda_1,\cdots, \lambda_{N+1}$.
It follows that $f_,\ldots,f_{N+1}$ are linearly dependent modulo
$Id(B^\sharp)$ and we have completed the proof.

\hfill $\Box$

\begin{lemma}\label{l8}Let $A=A(m,w)$ where $m\ge 2$ and $w$ is periodic or Sturmian word. Then
$$
l_n(A^\sharp) \le 4(m+1)(n+1)^{12}
$$
for all sufficiently large $n$.
\end{lemma}
{\em Proof}.
First note that the cocharacter of $A^\sharp$ lies in the strip of width $4$ that is, if $m_\lambda\ne 0$ in the decomposition
\begin{equation}\label{e5}
\chi_n(A^\sharp)=\sum_{\lambda\vdash n} m_\lambda\chi_\lambda
\end{equation}
then $h(\lambda)\le 4$. The number of partitions of $n$ of type $\lambda=(\lambda_1,\ldots,\lambda_k)$ with
$1\le k\le 4$ is less than $n^4$. By Lemma \ref{l6}
\begin{equation}\label{e6}
\dim W_n^{(4)}(A) \le 4(n+1) Comp_w(n).
\end{equation}
If $w$ is Sturmian word then $Comp_w(n)=n+1$. If $w$ is periodic then its complexity is finite and hence
$Comp_w(n)\le n+1$ for all sufficiently large $n$ in (\ref{e6}). In particular,
$$
\dim W_n^{(4)}(A) \le 4(m+1)(n+1)^2\le 4(m+1)n^2
$$
for all sufficiently large $n$. Applying Lemmas \ref{l5}, \ref{l6} and \ref{l7} we obtain
$$
m_\lambda \le\dim W_n^{(4)}(A^\sharp) \le 4(m+1)(n+1)^8
$$
for all $m_\lambda\ne 0$ in (\ref{e5}) and then
$$
l_n(A^\sharp) = \sum_{\lambda\vdash n}m_\lambda \le 4(m+1)n^4(n+1)^8 \le.
4(m+1)(n+1)^{12}.
$$
\hfill $\Box$

In the next step we shall find an upper bound for $\Phi(\lambda)$ provided that $m_\lambda\ne 0$ in the $n$th cocharacter
of $A^\sharp$.

\begin{lemma}\label{l9}
For any $\varepsilon>0$ there exists $n_0$ such that $m_\lambda=0$ in (\ref{e5}) if $n>n_0$ and
$$
\frac{\lambda_3}{\lambda_1}\ge \frac{\beta}{1-\beta} +\varepsilon
$$
where $\beta=\frac{1}{m+\alpha}$ and $\alpha$ is the slope of $w$.
\end{lemma}

{\em Proof}. First let $\lambda=(\lambda_1,\lambda_2,\lambda_3,1)\vdash n$. Inequality $m_\lambda\ne 0$ means that 
there exists a multilinear polynomial $g$ of degree $n$ depending on one alternating set of four variables,
$\lambda_3-1$ alternating sets of three variables and some extra variables and $g$ is not an identity of $A^\sharp$.
That is, there exists an evaluation $\varphi:F\{X\}^\sharp\to A^\sharp$ such that $\varphi(g)\ne 0$ and the set
$\{\varphi(x_1),\ldots,\varphi(x_n)\}$ contains at least $\lambda_3$ basis elements $b\in A$ and at most $\lambda_1$ elements 
$a\in A$. Obviously, $\varphi(g)=0$ if $\{\varphi(x_1),\ldots,\varphi(x_n)\}$ does not contain exactly one element 
$z_j^{(i)}\in A$.

Any non-zero product of basis elements of $A$ is the left-normed product of the type
$$
z_j^{(i)}a^{k_1} b^{l_1}\cdots a^{k_t}b^{l_t}
$$
where $k_1,\ldots,k_t,l_1,\ldots, l_t$ are equal to $0$ or $1$. More precisely, this product can be written in the form
\begin{equation}\label{e7}
z_j^{(i)} f(a,b)
\end{equation}
where
$$
f(a,b)=a^{t_0}ba^{t_1}b\cdots b a^{t_k}ba^{t_{k+1}}
$$
is an associative monomial on $a$ and $b$ and
$$
t_0=m+w_i-j, t_1=m+w_{i+1}-1,\ldots, t_k=m+w_{i+k}-1, t_{k+1} \le m+w_{i+k+1}-1.
$$
In particular, $\deg_bf=k+1$ and
$$
\deg_af=t_0+t_1+\cdots+t_{k+1}\ge t_1+\cdots+t_k=(m-1)k+ w_{i+1}+\cdots+w_{i+k}.
$$
The total degree of monomial (\ref{e7}) (i.e. the number of factors) is
$$
n=(m+1)k+ w_{i}+\cdots+w_{i+k}+t_{k+1}-j+1.
$$
Hence, $(m+1)k\ge n-(1+k)-m-1$ and $k\ge \frac{n-m-2}{m+2}$. In particular, $k$ grows with increasing  $n$.

It is known that
$$
\frac{w_{i+1}+\cdots+w_{i+k}}{k}\ge \alpha - \frac{C}{k}
$$
for some constant $C$ (see \cite{GMZAdvM}, Proposition 5.1 or \cite{Lotair}, Section 2.2). This implies that
$$
\deg_af>(m-1)k+k(\alpha-\delta)
$$
where $\delta=\frac{C}{k}$ and
$$
\frac{\deg_bf}{\deg_af} <\frac{1+\frac{1}{k}}{m-1+\alpha-\delta}.
$$
Since $\varphi(g)\ne 0$, at least one monomial of the type (\ref{e7}) in $\varphi(g)$ is non-zero. Therefore
\begin{equation}\label{e8}
\frac{\lambda_3}{\lambda_1}\le \frac{\deg_bf}{\deg_af}<\frac{1+\frac{1}{k}}{m-1+\alpha-\delta}.
\end{equation}
Since $\delta$ is an arbitrary small positive real number, one can choose $n_0$ such that
\begin{equation}\label{e9}
<\frac{1+\frac{1}{k}}{m-1+\alpha-\delta}<\frac{1}{m-1+\alpha}+\frac{\varepsilon}{2}
\end{equation}
for all $n\ge n_0$. Combining (\ref{e8}) and (\ref{e9}) we conclude that
\begin{equation}\label{e10}
\frac{\lambda_3}{\lambda_1}<\frac{1}{m-1+\alpha}+\frac{\varepsilon}{2}
\end{equation}
provided that $m_\lambda\ne 0$ in (\ref{e5}) and $n\ge n_0$. Note that $\frac{\beta}{1-\beta}=\frac{1}{m-1+\alpha}$,
hence we have completed the proof of our lemma in case $\lambda=(\lambda_1,\lambda_2,\lambda_3,1)$.

Slightly modifying previous arguments we get the proof of the inequality (\ref{e10}) for a partition
$\lambda=(\lambda_1,\lambda_2,\lambda_3)$ with three parts. The only difference is that non-identical
polynomial $g$ depends on at least $\lambda_3$ skewsymmetric sets of variables of order 3 but after evaluation, one 
of these variables can be replaced by $z_j^{(i)}$ and we get the inequality
$$
\frac{\lambda_3-1}{\lambda_1}\le \frac{\deg_bf}{\deg_af}
$$
instead of (\ref{e8}). Taking into account that $\lambda_1\to\infty$ if $n\to\infty$ we get the same conclusion and thus complete the proof.
\hfill $\Box$

For the lower bound of codimensions of $A^\sharp$ we need the following results.

Let $A=A(m,w)$ be an algebra defined by an integer $m\ge 2$ and by an infinite word $w=w_1w_2\ldots$ in the 
alphabet $\{0,1\}$. Then
\begin{equation}\label{e11}
z_1^{(1)}a^{i_1}ba^{i_2}b\cdots a^{i_r}b=z_1^{r+1}
\end{equation}
if $i_1=m-1+w_1, i_2=m-1+w_2,\ldots, i_r=m-1+w_r$. Otherwise the left hand side of (\ref{e11}) is zero.

\begin{lemma}\label{l11}
Let $\lambda=(j,\lambda_2,\lambda_3,1)$ be a partition of $n=j+mr+w_1+\cdots+w_r+1$ with 
$j\ge\lambda_2=(m-1)r+w_1+\cdots+w_r$, $\lambda_3=r$ or let $\lambda=(\lambda_1,j,\lambda_3,1)$   be
a partition of the same $n$ with $\lambda_1=(m-1)r+w_1+\cdots+w_r>j\ge \lambda_3=r$. Then $m_\lambda\ne 0$ in (\ref{e5}).
\end{lemma}

{\em Proof}. Recall that, given $S_n$-module $M$, the multiplicity of $\chi_\lambda$ in the character
$\chi(M)$ is non-zero if $e_{T_\lambda}M\ne 0$ for some Young tableaux $T_\lambda$ of shape $D_\lambda$.
The essential idempotent $e_{T_\lambda}\in FS_n$ is equal to
$$
e_{T_\lambda}=(\sum_{\sigma\in R_{T_\lambda}}\sigma)(\sum_{\tau\in C_{T_\lambda}}(sgn~\tau)\tau.
$$
Here $R_{T_\lambda}$ is the row stabilizer of $T_\lambda$, i.e. the subgroup of all $\sigma\in S_n$
permuting indices only inside rows of $T_\lambda$ and $C_{T_\lambda}$ is the column stabilizer of $T_\lambda$.
 
First let $\lambda_1=j\ge \lambda_2$. Denote $n_0=mr+w_1+\cdots+w_r+1$ and consider the Young tableaux $T_\lambda$
of the following type. Into the boxes of the 1st row of $D_\lambda$ we place $n_0+1,\ldots,n_0+j$ from left to right. 
Into the third row we insert $j_1=i_1+2,\ldots, j_r=i_1+\cdots+i_r+r+1$. (In fact,
$j_1,\ldots, j_r$ are the positions of $b$ in the product (\ref{e11})). 
Into the second row we insert from left to right
$j_1-1,\ldots,j_r-1, i_{r+1},\ldots, i_{\lambda_2}$ where $\{i_{r+1},\ldots, i_{\lambda_2}\}=
\{2,\ldots,n_0\}\setminus\{j_1-1,j_1,\ldots,j_r-1,j_r\}$ and into the unique box of the 4th row we put $1$.

Then
$$
e_{T_\lambda}(x_1,\ldots,x_n)=Sym_1Sym_2Sym_3Alt_1\cdots Alt_{\lambda_2}(x_1,\ldots,x_n)
$$
where
\begin{itemize}
\item[-]
$Alt_1$ is the alternation on $\{1,j_1-1,j_1,n_0+1\}$;

\item[-]
$Alt_k$ is the alternation on $\{j_k-1,j_k,n_0+k\}$ if $2\le k \le r$;

\item[-]
$Alt_k$ is the alternation on $\{i_k,n_0+k\}$ if $r<k\le \lambda_2$;

\item[-]
$Sym_1$ is the symmetrization on $\{n_0+1,\ldots, n_0+j\}$;

\item[-]
$Sym_2$ is the symmetrization on $\{j_1,\ldots,j_r\}$;

\item[-]
$Sym_3$ is the symmetrization on $\{2,\ldots,n\}\setminus \{j_1,\ldots,j_r\}$.
\end{itemize}

After an evaluation
$$
\varphi(x_1)=z_1^{(1)}, 
\varphi(x_{n_0+1})=\cdots=\varphi(x_{n_0+j})=1\in A^\sharp, 
\varphi(x_{j_1})=\cdots=\varphi(x_{j_r})=b
$$
and
$$
\varphi(x_i)=a~if~ i\ne 1,j_1,\ldots,j_r,n_0+1,\ldots,n_0+j
$$
we have
$$
\varphi(e_{T_\lambda}(x_1\cdots x_n))= j! r! (n_0-r-1)! z_1^{(r+1)}\ne 0,
$$
hence $m_\lambda\ne 0$ in (\ref{e5}).

Similarly, filling up the second row of $T_\lambda$ by $n_0+1,\ldots, n_0+j$ in case
$\lambda_1=(m-1)+w_1+\cdots+w_r> j\ge \lambda_3=r$ we prove that $e_{T_\lambda}(x_1\cdots x_n)$
is not an identity of $A^\sharp$.
\hfill $\Box$

Recall that, given $0\le\beta\le 1$,
$$
\Phi_0(\beta)=\Phi(\beta,1-\beta)=\frac{1}{\beta^\beta(1-\beta)^{1-\beta}}.
$$

\begin{lemma}\label{l12}
Let $A=A(m,w)$ be an algebra defined for an integer $m\ge 2$ and a Sturmian or periodic word $w$
with  slope $\alpha$. Let also $\beta=\frac{1}{m+\alpha}$. Then for any $\varepsilon>0$ there exist 
a constant $C$, positive integers $n_1< n_2\ldots$ and partitions $\lambda^{(i)}\vdash n_i$ such that for some large enough $i_0$ the following properties hold:
\begin{itemize}
\item[1)]
$\vert\Phi(\lambda^{(i)})-\Phi_0(\beta)-1\vert < \varepsilon\quad {\rm for~ all}\quad i\ge i_0$;

\item[2)]
$n_{i+1}-n_i < C$ for all $i\ge i_0$;

\item[3)]
$m_\lambda^{(i)}\ne 0$ in $\chi _{n_i}(A^\sharp)$ for all $i\ge i_0$.
\end{itemize}
\end{lemma}

{\em Proof}.  Note that $\beta<\frac{1}{2}$ since $\alpha >0$.
First take an arbitrary $r\ge 1$, $n=mr+w_1+\cdots+w_r$ and
$\lambda=(\lambda_1,\lambda_2)$, where $\lambda_1=(m-1)r+w_1+\cdots+ w_r, \lambda_2=r$. We set
$$
x_1=\frac{\lambda_1}{n}=\frac{m-1+\frac{w_1+\cdots+ w_r}{r}}{m+\frac{w_1+\cdots+ w_r}{r}},
$$
$$
x_2=\frac{\lambda_2}{n}=\frac{1}{m+\frac{w_1+\cdots+ w_r}{r}}.
$$
As it was mentioned in the proof of Lemma \ref{l9} (see also \cite{GMZAdvM}, Proposition 5.1 or \cite{Lotair},
Section 2.2) there exists a constant $C_1$ such that
\begin{equation}\label{e11a}
\vert\frac{w_1+\cdots+ w_r}{r}-\alpha\vert <\frac{C_1}{r}.
\end{equation}
Hence for any $\varepsilon_1>0$  we can find $r_0$ such that
\begin{equation}\label{e12}
\vert\Phi(\lambda)-\Phi_0(\beta)\vert < \varepsilon_1
\end{equation}
for all $r\ge r_0$ since $\Phi(z_1,z_2)$ is a continuous function and $(x_1,x_2)\to (1-\beta,\beta)$
when $r\to\infty$.

Now using Lemma \ref{l3} and Lemma \ref{l4}, given $\varepsilon_2>0$, we insert one extra 
row into $D_\lambda$ that is, we construct  a partition $\mu=(\mu_1,\mu_2,\mu_3)$ of $n_0=nk$
such that
\begin{equation}\label{e14}
\vert\Phi(\lambda)-\Phi(\mu)-1\vert < \varepsilon_2.
\end{equation}

We have three options. Either $\mu_1$ is a new row (that is, $(\mu_2,\mu_3)=
(q\lambda_1,q\lambda_2)$ or $\mu_2$ is a new row (that is, $(\mu_1,\mu_3)=
(q\lambda_1,q\lambda_2)$ or $\mu_3$ is a new row (that is, $(\mu_1,\mu_2)=
(q\lambda_1,q\lambda_2)$.

First we exclude the third case.
Suppose that $(\mu_1,\mu_2)=(q\lambda_1,q\lambda_2)$.
Recall that by Lemma \ref{l3} the maximal value of $\Phi(tz_1,tz_2,1-t)$ is achieved if 
$$
t=\frac{\Phi(z_1,z_2)}{1+\Phi(z_1,z_2)}.
$$
Since $\Phi(z_1,z_2)<2$ if $\beta<\frac{1}{2}$ we obtain that $1-t> \frac{1}{3}$. For Lemma \ref{l4} this means
that the new row of $D_\mu$ cannot be the third row that is, case $(\mu_1,\mu_2)=(q\lambda_1,q\lambda_2)$
is impossible.

Now let $(\mu_2,\mu_3)=(q\lambda_1,q\lambda_2)$. We exchange $\mu$ to $\mu'$ in the following way. We set
$\mu_2'=qr(m-1)+w_1+\cdots+w_{qr}$ and take $\mu'=(\mu_1,\mu_2',\mu_3)$. Then $\mu'\vdash n'$ where
$$
n'-n_0=\mu_2'-\mu_2=w_1+\cdots+w_{qr}-q(w_1+\cdots+w_{r}).
$$
Using again inequality (\ref{e11a}) we get
\begin{equation}\label{e15}
|n'-n_0|<  C_1(q+1).
\end{equation}
Inequality (\ref{e15}) also shows that $\mu_1\ge\mu_2'\ge\mu_3$ if $n$ is sufficiently large
and our construction of partition $\mu$ is correct.

Clearly, $\vert \Phi(\mu)-\Phi(\mu') \vert\to 0$ if $n\to\infty$ and
\begin{equation}\label{e16}
\vert \Phi(\mu)-\Phi(\mu') \vert <\varepsilon_3
\end{equation}
for any fixed $\varepsilon_3 >0$ for all sufficiently large $r$ (and $n$). Starting from this sufficiently large
$r$ we denote $n_r=n'+1$ and take $\lambda^{(r)}\vdash n_r$, $\lambda^{(r)}=(\mu_1,\mu',\mu_3,1)$. All
preceding $n_1,\ldots,n_{r-1}$ and $\lambda^{(1)},\ldots,\lambda^{(r-1)}$ we choose 
in an arbitrary way.

Since $\mu_3=qr$, by Lemma \ref{l11} the multiplicity of the irreducible character $\lambda^{(r)}$ in
$\chi_{n_r}(A^\sharp)$ is not equal to zero and $|n_r-kn|< C_2=C_1(q+1)+1$ by (\ref{e15}) since $n_0=nk$.
It is not difficult to see that in this case
\begin{equation}\label{e17}
\vert \Phi(\mu')-\Phi(\lambda^{(r)}) \vert <\varepsilon_4
\end{equation}
for any fixed $\varepsilon_4 >0$ if $r$ (and the corresponding $n$) is sufficiently large. Combining
(\ref{e12}), (\ref{e14}), (\ref{e16}) and (\ref{e17}) we see that $\lambda^{(r)}$ satisfies
conditions 1) and 3) of our lemma. Finally, consider the difference between $n_r$ and $n_{r+1}$
provided that all $n_{r+1}, n_{r+2},\ldots$ are constructed by the same procedure. 
That is, we take
$$
\bar n = m(r+1)+w_1+\cdots+w_{r+1}+1
$$
and obtain $n_{r+1}$ satisfying the same condition
$$
\vert n_{r+1}- k\bar n \vert < C_2.
$$
On the other hand, $\bar n -kn=k(m+w_{r+1}) \le k(m+1)$ and 
$\vert kn-n_r\vert < C_2$. Hence we have
$$
\vert n_{r+1} - n_r \vert < C=2C_2+k(m+1).
$$
This latter inequality completes the proof of our lemma if $(\mu_2,\mu_3)=(q\lambda_1,q\lambda_2)$.
Arguments in the case $(\mu_1,\mu_3)=(q\lambda_1,q\lambda_2)$ are the same.
\hfill $\Box$

\section{The main result}

Now we are ready to prove the main result of the paper.

\begin{theorem}\label{t1}
Let $w=w_1w_2\ldots$ be Sturmian or periodic word and let $A=A(m,w)$, $m\ge 2,$ be an algebra defined 
by $m$ and $w$ in (\ref{e01}) - (\ref{e04}). If $A^\sharp$ is the algebra obtained from $A$ 
by adojining an external unit then PI-exponent of $A^\sharp$ exists and
$$
exp(A^\sharp)=1+exp(A).
$$
\end{theorem}

{\em Proof.}
Let $\alpha=\pi(w)$ be the slope of $w$ and let $\beta=\frac{1}{m+\alpha}$. Recall that $exp(A)=\Phi_0(\beta)$
where
$$
\Phi_0(\beta)=\frac{1}{\beta^\beta(1-\beta)^{1-\beta}}
$$
(\cite{GMZAdvM}). First we prove that for any $\delta>0$ there exists $N$ such that
\begin{equation}\label{e18}
\Phi(\lambda)<\Phi_0(\beta)+1+\delta
\end{equation}
as soon as $\lambda$ is a partition of $n\ge N$ with $m_\lambda\ne 0$ in $\chi_n(A^\sharp)$.

By Lemma \ref{l9}, for any $\varepsilon>0$ there exists $n_0$ such that
\begin{equation}\label{e23a}
\frac{\lambda_3}{\lambda_1}<\frac{\beta}{1-\beta}+\varepsilon
\end{equation}
if $n\ge n_0, \lambda\vdash n$ and $m_\lambda\ne 0$. If $\lambda=(n)$ or $\lambda=(\lambda_1,\lambda_2)$
then by the hook formula for dimensions of irreducible $S_n$-representations it follows that 
$\deg \chi_\lambda\le 2^n$. Then by Lemma \ref{l1}
$$
\Phi(\lambda)\le 2\sqrt[n]{n^6}
$$
and (\ref{e18}) holds for all sufficiently large $n$ since $1\le\Phi_0(\beta)\le 2$.

Let $\lambda=(\lambda_1,\lambda_2,\lambda_3)$. Denote $\mu=(\lambda_1,\lambda_3)\vdash n'$, where
$n'=n-\lambda_2$. If $x_1=\frac{\lambda_1}{n'}, x_2=\frac{\lambda_3}{n'}$ then
$$
\Phi(\mu)=\Phi(x_1,x_2)=\Phi_0(x_2)
$$
and
$$
x_2\le \varphi(\varepsilon)=\frac{\beta+(1-\beta)\varepsilon}{1+(1-\beta)\varepsilon}
$$
as follows from (\ref{e23a}). Since
$$
\lim_{n\to\infty}\varphi(\varepsilon)=\beta
$$
and $\Phi_0$ is continuous, there exist $N$ and $\varepsilon$ such that 
$\Phi(\mu)<\Phi_0(\beta)+\delta$ for all $n>N$.
Then by Lemma \ref{l3}
$$
\Phi(\lambda)\le \Phi(\mu)+1< \Phi_0(\beta)+1+\delta.
$$

Now consider the case $\lambda=(\lambda_1,\lambda_2,\lambda_3,1)$. Excluding
the second row of diagram $D_\lambda$
we get a partition $\mu=(\mu_1,\mu_2,1)=(\lambda_1,\lambda_3,1)$ of $n'=n-\lambda_2$ with
$$
\frac{\mu_2}{\mu_1}<\frac{\beta}{1-\beta}+\varepsilon.
$$

Consider also partition $\mu'=(\mu_1,\mu_2)$ of $n'-1$. As before, given $\delta>0$, one can find $n_0$
such that
$$
\Phi(\mu')<\Phi_0(\beta)+\frac{\delta}{2}
$$
provided that $n\ge n_0$.

Since $\Phi$ is continuous, for all sufficiently large $n$ (and $n'$) we have
$$
\Phi(\mu) < \Phi_0(\beta)+\delta.
$$
Applying again Lemma \ref{l3} we get (\ref{e18}). It now follows from (\ref{e18}), Lemmas \ref{l1} and  \ref{l8}  that
$$
c_n(A^\sharp)=\sum_{\lambda\vdash n}m_\lambda\deg\chi_\lambda \le(\Phi_0(\beta)+1+\delta)^n l_n(A^\sharp)
\le 4(m+1)(n+1)^{13} (\Phi_0(\beta)+1+\delta)^n.
$$
Hence
$$
\overline{exp}(A^\sharp)=\limsup_{n\to\infty}\sqrt[n]{c_n(A^\sharp)}\le\Phi_0(\beta)+\delta+1
$$
for any $\delta>0$ that is,
\begin{equation}\label{e19}
\overline{exp}(A^\sharp)\le \Phi_0(\beta)+1=exp(A)+1.
\end{equation}

Now we find a lower bound for codimensions. Since
$$
c_n(A^\sharp)\ge\deg\chi_\lambda\ge \frac{\Phi(\lambda)^{n}}{n^{20}}
$$
by Lemma \ref{l1} if $m_\lambda\ne 0$ in $\chi_n(A^\sharp)$, then by Lemma \ref{l12} for any $\varepsilon>0$ there exists a sequence $n_1<n_2<\ldots$ such that
$$
c_{n_i}(A^\sharp)\ge \frac{1}{n_i^{20}}(\Phi_0(\beta)+1-\varepsilon)^{n_i},~i=1,2,\ldots
$$
and $n_{i+1}-n_i<C=const$ for all $i\ge 1$. Note that the sequence $\{c_n(R)\}$ is non-decreasing for
any unital algebra $R$. Then
\begin{equation}\label{e20}
\underline{exp}(A^\sharp)=\liminf_{n\to\infty}\sqrt[n]{c_n(A^\sharp)}\ge \Phi_0(\beta)+1.
\end{equation}
Now (\ref{e19}) and (\ref{e20}) complete the proof of the theorem.
\hfill $\Box$

\begin{corollary}\label{c1}
For any real numbers $\gamma\in [2,3]$ there exists an algebra $A$ with 1 such that $exp(A)=\gamma$.
\end{corollary}
\hfill $\Box$

As it was mentioned in the preliminaries, PI-exponents of finite dimensional algebras form
a dense subset in the interval $[1,2]$. Hence we get the following

\begin{corollary}\label{c2}
For any real numbers $\beta<\gamma\in [2,3]$ there exists a finite dimensional algebra $B$ 
with 1 such that $\beta\le exp(B)\le\gamma$. In particular, PI-exponents of finite
dimensional unital algebras form a dense subset in the interval $[2,3]$.
\end{corollary}
\hfill $\Box$

\end{document}